\documentclass[a4paper]{jpconf}
\usepackage{graphicx}
\usepackage{amssymb}
\usepackage{amsmath}
\usepackage{epsf}
\usepackage{bm}
\usepackage{cite}

\newtheorem{theorem}{Theorem}

\begin{document}

\title{Symmetries of second-order PDEs and conformal Killing vectors}
\author{Michael Tsamparlis$^1$ and Andronikos Paliathanasis$^{2,3}$}

\begin{abstract}
We study the Lie point symmetries of a general class of partial differential
equations (PDE) of second order. An equation from this class naturally
defines a second-order symmetric tensor (metric). In the case the PDE is
linear on the first derivatives we show that the Lie point symmetries are
given by the conformal algebra of the metric modulo a constraint involving
the linear part of the PDE. Important elements in this class are the
Klein--Gordon equation and the Laplace equation. We apply the general
results and determine the Lie point symmetries of these equations in various
general classes of Riemannian spaces. Finally we study the type II\ hidden
symmetries of the wave equation in a Riemannian space with a Lorenzian
metric.
\end{abstract}

\address{$^1$ Department of Physics, Section of Astronomy, Astrophysics and
Mechanics, University of Athens, Panepistemiopolis, Athens 157 83, Greece}
\address{$^2$ Dipartimento di Fisica, Universit\`{a} di Napoli ``Federico
II'' Complesso universitario Monte S.~Angelo Via, Cintia - 80126 Napoli, Italy}
\address{$^3$ INFN - Istituto Nazionale di Fisica Nucleare sezione di Napoli Complesso universitario Monte S.~Angelo Via Cintia - 80126 Napoli,
Italy} \ead{mtsampa@phys.uoa.gr, paliathanasis@na.infn.it}

\section{Introduction}

In theoretical physics one has two main tools to study the properties of
evolution of dynamical systems (a) Symmetries of the equations of motion and
(b) Collineations (symmetries) of the background space, where evolution
takes place. It is well known that both these tools have the following
common characteristics:

1) they form a Lie algebra;

2) they do not fix uniquely either the dynamical system or the space.

The natural question to be to asked is if these two algebras are related and
in what way. Equivalently, one may state the question as follows:

\emph{To what degree and how the space modulates the evolution of dynamical
systems in it? That is, a dynamical system is free to evolve at will in a
given space or it is constrained to do so by the very symmetry structure of
the space?}\

This question has been answered many years ago by the Theory of Relativity
with the Equivalence Principle, that is, the requirement that free motion in
a given gravitational field occurs along the geodesics of the space. However
as obvious as this point of view may appear to be it is not easy to
comprehend and accept! So let us give a precise formulation now\footnote{%
This point has been raised during an illuminating discussion with Prof P G
Leach in 2000 in Athens while we were driving to the Poseidon Temple in Cape
Sounion.}.

In a Riemannian space the affinely parameterized geodesics are determined
uniquely by the metric. The geodesics are a set of homogeneous ordinary
differential equations (ODE) linear in the highest order term and
quadratically non-linear in the first order terms. A system of such ODEs is
characterized (not fully) by its Lie point symmetries. On the other hand a
metric is characterized (again not fully) by its collineations. Therefore it
is reasonable one to expect that the Lie point symmetries of the system of
geodesic equations of a metric will be closely related with the
collineations of the metric. That such a relation exists it is easy to see
by the following simple example. Consider on the Euclidian plane a family of
straight lines parallel to the $x$-axis. These curves can be considered
either as the integral curves of the ODE\ $\frac{d^{2}y}{dx^{2}}=0$ or as
the geodesics of the Euclidian metric $dx^{2}+dy^{2}$. Subsequently consider
a symmetry operation defined by a reshuffling of these lines without
preserving necessarily their parametrization. According to the first
interpretation this symmetry operation is a Lie symmetry of the ODE $\frac{%
d^{2}y}{dx^{2}}=0$ and according to the second interpretation it is a
(special) projective symmetry of the Euclidian two-dimensional space.

What has been said for a Riemannian space can be generalized to an affine
space in which there is only a linear connection. In this case the geodesics
are called autoparallels (or paths) and they comprise again a system of ODEs
linear in the highest order term and quadratically non-linear in the first
order terms. In this case one is looking for relations between the Lie point
symmetries of the autoparallels and the projective collineations of the
connection.

A Lie point symmetry of an ordinary differential equation (ODE) is a point
transformation in the space of variables which preserves the set of
solutions of the ODE \cite{Ibrag,Stephani,Bluman}. If we look at these
solutions as curves in the space of variables, then we may equivalently
consider a Lie point symmetry as a point transformation which preserves the
set of the solution curves. Applying this observation to the geodesic curves
in a Riemannian (affine) space, we infer that the Lie point symmetries of
the geodesic equations in any Riemannian space are the automorphisms which
preserve the set of these curves. However we know from Differential Geometry
that the point transformations of a Riemannian (affine)\ space which
preserve the set of geodesics are the projective transformations. Therefore
it is reasonable to expect a correspondence between the Lie point symmetries
of the geodesic equations and the projective algebra of the metric defining
the geodesics.

The equation of geodesics in an arbitrary coordinate frame is a second-order
ODE of the form%
\begin{equation}
\ddot{x}^{i}+\Gamma _{jk}^{i}\dot{x}^{j}\dot{x}^{k}+F^{i}(x,\dot{x})=0,
\label{de.0}
\end{equation}%
where $F^{i}(x,\dot{x})$ are arbitrary functions of their arguments and the
functions $\Gamma _{jk}^{i}$ are the connection coefficients of the space.
Equivalently equation (\ref{de.0}) is also the equation of motion of a
dynamical system moving in a Riemannian (affine) space under the action of a
velocity dependent force. According to the above argument we expect that the
Lie point symmetries of the ODE (\ref{de.0}) for given functions $F^{i}(x,%
\dot{x})$ will be related with the collineations of the metric. As it will
be shown in this case the Lie symmetries of (\ref{de.0}) determine a
subalgebra of the special projective algebra of the space. The specific
subalgebra is selected by means of certain constraint conditions involving
geometric quantities of the space and the function $F(x^{i},\dot{x}^{j})$
\cite%
{PrinceCrampin(1984)1,Aminova2006,Aminova2010,ND2010,TsamparlisGRG,TsamparlisGRG2}%
.

The determination of the Lie point symmetries of a given system of ODEs
consists of two steps: (a)\ determination of the conditions which the
components of the Lie symmetry vector must satisfy and (b) solution of the
system of these conditions. Step (a) is formal and it is outlined, e.g., in
\cite{Ibrag,Stephani,Bluman}. The second step is the key one and, for
example, in higher dimensions, where one has a large number of equations,
the solution can be quite involved and perhaps impossible by algebraic
computing. However, if one expresses the system of Lie symmetry conditions
of (\ref{de.0}) in terms of collineation (i.e. symmetry) conditions of the
metric, then the determination of Lie point symmetries is transferred to the
geometric problem of determining the special projective group of the metric
\cite{TsamparlisGRG2}. In this field there is a significant amount of
knowledge from Differential Geometry waiting to be used. Indeed the
projective symmetries are already known for many spaces or they can be
determined by existing general theorems. For example the projective algebra
and all its subalgebras are known for the important case of spaces of
constant curvature \cite{Barnes1993} and in particular for the flat spaces.
This implies that, the Lie symmetries of the Newtonian dynamical systems as
well as those of Special Relativity can be determined by simple
differentiation from the known projective algebra of these spaces!

What has been said for the Lie point symmetries of (\ref{de.0}) applies also
to Noether point symmetries (provided (\ref{de.0}) follows from a
Lagrangian). The Noether point symmetries are Lie point symmetries which
satisfy the additional constraint%
\begin{equation}
X^{\left[ 1\right] }L+L\frac{d\xi }{dt}=\frac{df}{dt}.  \label{L2p.3}
\end{equation}%
The Noether point symmetries form a subalgebra of the Lie point symmetry
algebra. In accordance with the above this implies that the Noether point
symmetries will be related with a subalgebra of the special projection
algebra of the space where `motion' occurs. As it has been shown this
subalgebra is the homothetic algebra of the space \cite{TsamparlisGRG2}. It
is well known that to each Noether point symmetry it is associated a
conserved current (i.e. a Noether first integral). This leads to the
important conclusion that the (standard) conserved quantities of a dynamical
system depend \textit{on the space it moves and the type of force }$F(x^{i})$
\textit{\ which modulates the motion}. In particular in `free fall', that is
when $F(x^{i})=0,$ the orbits are affinely parametrized geodesics and the
\emph{geometry} of the space is the sole factor which determines the
conserved quantities of motion. This conclusion is by no means trivial and
means that the space where motion occurs is not a simple carrier of the
motion but it is the major modulator of the evolution of a dynamical system.
In other words there is a strong and deep relation between Geometry of the
space and Physics (motion) in that space!

The above scenario can be generalized to partial differential equations
(PDEs). Obviously, in this case a global answer is not possible. However, it
can be shown that for many interesting PDEs the Lie point symmetries are
indeed obtained from the collineations of the metric. Pioneering work in
this direction is the work of Ibragimov \cite{Ibrag}. Recently, Bozhkov et
al.~\cite{Boskov} studied the Lie and the point Noether symmetries of the
Poisson equation and have shown that the Lie symmetries of the Poisson PDE
are generated from the conformal algebra of the metric. This result can be
generalized and it has been shown \cite{JGP1} that for a general class of
PDEs of second order in an $n$-dimensional Riemannian space, there is a
close relation between the Lie point symmetries and the conformal algebra of
the space. Examples of such PDEs include some important equations as: the
heat equation, the Klein--Gordon equation, the Laplace equation, the Schr%
\"{o}dinger equation and others \cite{JGP1,JGP2,AnIJGMMP,MTIJGMMP}. In what
follows we discuss in a rather systematic way the aforementioned ideas. The
plan of the paper is as follows.

In section~\ref{Prel} we give the basic definitions and properties of the
collineations of space times and the Lie symmetries of DEs. In section~\ref%
{ScPDE} we study the Lie symmetries of a generic family of second-order
partial differential equations and we prove that when the second-order
partial differential equation is linear on the derivatives, the
corresponding Lie symmetries are generated by the conformal algebra of the
underlying geometry.\ Furthermore, in sections~\ref{PoissonSym} and~\ref%
{PoissonSym1} we apply these general result in order to determine the
general form of the Lie symmetry vector of the Poisson equation, of the
Klein--Gordon equation and of the Laplace equation. In section~\ref%
{Reduction} we study the application of the conformal Killing vectors in the
Laplace equation is some special Riemannian spaces and we study the origin
of the type II hidden symmetries. Finally in section~\ref{Application} we
apply the previous results in the case of the Laplace equation in the 1+3
wave equation and in Bianchi~I spacetimes.

\section{Preliminaries}

\label{Prel}

In this section we give the basic definitions and properties of the
collineations of spacetimes and of the point symmetries of differential
equations.

\subsection{Collineations of Riemannian spaces}

A collineation in a Riemannian space of dimension $n$ is a vector field $%
\mathbf{X}$ which satisfies an equation of the form~$\mathcal{L}_{X}\mathbf{A%
}=\mathbf{B}$, where $\mathcal{L}_{X}$ is the Lie derivative with respect to
the vector field~$\mathbf{X}$, $\mathbf{A}$ is a geometric object (not
necessary\ a tensor field) defined in terms of the metric and its
derivatives, and $\mathbf{B}$ is an arbitrary tensor field with the same
tensor indices as the geometric object $\mathbf{A}.$ The classification of
the collineations of Riemannian manifolds can be found in \cite{Katzin}. In
the following we are interested in the collineations of the metric tensor,
i.e.~$\mathbf{A=}g_{ij}~$of the Riemannian space.

A vector field $\mathbf{X}$ is a conformal Killing vector (CKV) of $g_{ij}$
if the following condition holds\footnote{%
In the following we use the Einstein notation.}:%
\begin{equation*}
\mathcal{L}_{\mathbf{X}}g_{ij}=2\psi \big(x^{k}\big) g_{ij},
\end{equation*}%
where $\psi \big(x^{k}\big) =\frac{1}{n}\mathbf{X}_{;i}^{i}$. In case $\psi
_{;ij}=0$, $\mathbf{X}$ is called special CKV (sp.CKV), if $\psi \big(
x^{k}\big) =$~constant, the vector field $X$ is called homothetic (HV) and
if $\psi \big(x^{k}\big) =0,$ the field $\mathbf{X}$ is called a Killing
vector (KV).

The CKVs of the metric $g_{ij}$ form a Lie algebra which is called the
conformal algebra ($CA)$ of the metric $g_{ij}~$(CKVs). The conformal
algebra contains two subalgebras, the Homothetic algebra $(HA)$~and the
Killing algebra $\left( KA\right) $. These algebras are related as $%
KA\subseteq HA\subseteq CA.$

Two metrics $g_{ij}$ and $\bar{g}_{ij}$ are conformally related if there
exists a function $N^{2}\big(x^{k}\big) $ such as $\bar{g}_{ij}=N^{2}\big(
x^{k}\big) g_{ij}.$ If $\mathbf{X}$ is a CKV\ of the metric $\bar{g}_{ij}$
so that $L_{\mathbf{X}}\bar{g}_{ij}=2\bar{\psi}\bar{g}_{ij},$ then $\mathbf{X%
}$ is also a CKV of the metric $g_{ij}$, that is $L_{\mathbf{X}}g_{ij}=2\psi
g_{ij}$ with conformal factor $\psi \big(x^{k}\big) $; the two conformal
factors are related as follows
\begin{equation*}
\psi =\bar{\psi}N^{2}-NN_{,i}^{i}\mathbf{X}^{,i}.
\end{equation*}%
The last relation implies that two conformally related metrics have the same
conformal algebra, but with different subalgebras; that is, a KV for one may
be proper CKV for the other.

A special class of conformally related spaces are the conformally flat
spaces. A space $V^{n}$ is conformally flat if the metric $g_{ij}$ of $V^{n}$
satisfies the relation $g_{ij}=N^{2}s_{ij}$, where $s_{ij}\,$~is the metric
of a flat space which has the same signature as $g_{ij}$. The maximal
dimension of the conformal algebra of a $n$-dimensional metric $\left(
n>2\right) $ is $\frac{1}{2}\left( n+1\right) \left( n+2\right) $ and in
that case the space is conformally flat. Moreover, if the conformally flat
space $V^{n}$ admits a $\frac{1}{2}n\left( n+1\right) $-dimensional Killing
algebra then $V^{n}$ is a space of constant curvature and admits a proper HV
if and only if the space is flat.

Furthermore, if for a CKV $\mathbf{X}$ of the Riemannian space $V_{G}^{n}$
the condition $\mathbf{X}_{\left[ i;j\right] }=0$ holds, i.e.~$\mathbf{X}%
_{i;j}=\psi \big(
x^{k}\big) g_{ij}$, then the CKV will be called gradient CKV. In this case
there exists a~coordinate system in which the line element of the metric
which defines the Riemannian space~$V_{G}^{n}$ is%
\begin{equation*}
ds^{2}=dx^{n}+f^{2}\big( x^{n}\big) h^{AB}\big( x^{A}\big) dx^{A}dx^{B},
\end{equation*}%
where $A,B=1,2,\dots,n-1$ \cite{Daftardar}. In these coordinates the
gradient CKV is $\mathbf{X}=f\left( x^{n}\right) \partial _{x^{n}}$ with
conformal factor $\psi =f_{,x^{n}}$. In the case when $f\left( x^{n}\right)
=x^{n},~$we have $\psi =1$; hence, $\mathbf{X}$ becomes gradient HV and if $%
f\left( x^{n}\right) =f_{0}$, $\mathbf{X}$ becomes gradient KV.

\subsection{Point symmetries of differential equations}

\label{Preliminaries}

A partial differential equation (PDE) is defined by a function $%
H=H(x^{i},u^{A},u_{,i}^{A},u_{,ij}^{A})$ in the jet space $\bar{B}_{\bar{M}}$%
, where $x^{i}$ are the independent variables and $u^{A}$ are the dependent
variables. The infinitesimal point transformation
\begin{align}
\bar{x}^{i}& =x^{i}+\varepsilon \xi ^{i}\big(x^{k},u^{B}\big),  \label{pr.01}
\\
\bar{u}^{A}& =\bar{u}^{A}+\varepsilon \eta ^{A}\big(x^{k},u^{B}\big),
\label{pr.02}
\end{align}%
has the infinitesimal symmetry generator
\begin{equation}
\mathbf{X}=\xi ^{i}\big(x^{k},u^{B}\big)\partial _{x^{i}}+\eta ^{A}\big(%
x^{k},u^{B}\big)\partial _{u^{A}}.  \label{pr.03}
\end{equation}

The generator $\mathbf{X}$ of the infinitesimal transformation (\ref{pr.01}%
),~(\ref{pr.02}) is called a Lie point symmetry of the PDE $H$ if there
exists a function $\lambda $ such that the following condition holds \cite%
{Ibrag,Stephani}
\begin{equation}
\mathbf{X}^{[n]}(H)=\lambda H,~~\bmod H=0,  \label{pr.04}
\end{equation}%
where
\begin{equation}
\mathbf{X}^{[n]}=\mathbf{X}+\eta _{i}^{A}\partial _{\dot{x}^{i}}+\eta
_{ij}^{A}\partial _{u_{ij\dots i_{n}}^{A}}+\dots+\eta _{i_{1}i_{2}\dots
i_{n}}^{A}\partial _{u_{i_{1}i_{2}\dots i_{n}}^{A}}  \label{pr.05}
\end{equation}%
is the $n$-th prolongation of $\mathbf{X}$ and
\begin{equation}
\eta _{i}^{A}=\eta _{,i}^{A}+u_{,i}^{B}\eta _{,B}^{A}-\xi
_{,i}^{j}u_{,j}^{A}-u_{,i}^{A}u_{,j}^{B}\xi _{,B}^{j}  \label{pr.06}
\end{equation}%
with

\begin{align}
\eta _{ij}^{A}& =\eta _{,ij}^{A}+2\eta _{,B(i}^{A}u_{,j)}^{B}-\xi
_{,ij}^{k}u_{,k}^{A}+\eta _{,BC}^{A}u_{,i}^{B}u_{,j}^{C}-2\xi
_{,(i|B|}^{k}u_{j)}^{B}u_{,k}^{A}  \notag \\
& -\xi _{,BC}^{k}u_{,i}^{B}u_{,j}^{C}u_{,k}^{A}+\eta
_{,B}^{A}u_{,ij}^{B}-2\xi _{,(j}^{k}u_{,i)k}^{A}-\xi _{,B}^{k}\left(
u_{,k}^{A}u_{,ij}^{B}+2u_{(,j}^{B}u_{,i)k}^{A}\right).  \label{pr.07}
\end{align}

Lie symmetries of differential equations can be used in order to determine
invariant solutions or transform solutions into solutions \cite{Bluman}.
From condition (\ref{pr.04}) one defines the Lagrange system%
\begin{equation}
\frac{dx^{i}}{\xi ^{i}}=\frac{du}{\eta }=\frac{du_{i}}{\eta _{\left[ i\right]
}}=\dots=\frac{du_{ij\dots i_{n}}}{\eta _{\left[ ij\dots i_{n}\right] }}
\label{pr.08}
\end{equation}%
whose solution provides the characteristic functions
\begin{equation*}
W^{[ 0] }\big( x^{k},u\big),\quad W^{[ 1] i}\big( x^{k},u,u_{,i}\big),\quad
\dots, \quad W^{[ n] }\big( x^{k},u,u_{,i},\dots,u_{,ij\dots i_{n}}\big).
\end{equation*}%
The solution $W^{\left[ n\right] }~$of the Lagrange system (\ref{pr.08}) is
called the $n$-th order invariant of the Lie symmetry vector (\ref{pr.03})
and holds $X^{\left[ n\right] }W^{\left[ n\right] }=0$.

The application of a Lie symmetry to a PDE $H$ leads to a new differential
equation $\bar{H}$ which is different from $H$ and it is possible that it
admits Lie symmetries which are not Lie symmetries of~$H$. These Lie point
symmetries are called Type~II hidden symmetries. It has been shown in \cite%
{Govinger} that if $X_{1},X_{2}$ are Lie point symmetries of the original
PDE with commutator $\left[ X_{1},X_{2}\right] =cX_{1}~$where $c$ is a
constant, then reduction by $X_{2}$ results in $X_{1}$ being a point
symmetry of the reduced PDE $\bar{H}$, while reduction by $X_{1}$ results in
a PDE $\bar{H}$ for which $X_{2}$ is not a Lie point symmetry.

In the following section we study the Lie point symmetries of a general type
of second-order PDEs.

\section{Lie symmetries of second-order PDEs and CKVs}

\label{ScPDE}

It is interesting to examine if the close relation of the Lie and the
Noether point symmetries of the second-order ODEs of the form (\ref{de.0})
with the collineations of the metric is possible to be carried over to
second-order partial differential equations (PDEs). Obviously it will not be
possible to give a complete answer, due to the complexity of the study and
the great variety of PDEs.

We consider the second-order PDEs of the form
\begin{equation}
A^{ij}u_{ij}-F(x^{i},u,u_{i})=0  \label{GPE.0}
\end{equation}%
for which at least one of the $A^{ij}$ is nonzero and derive the point Lie
symmetry conditions.

The symmetry condition (\ref{pr.04}) when applied to (\ref{GPE.0}) gives
\begin{equation}
A^{ij}\eta _{ij}^{\left( 2\right) }+\left( XA^{ij}\right)
u_{ij}-X^{[1]}(F)=\lambda (A^{ij}u_{ij}-F)  \label{GPE.13}
\end{equation}
that leads to
\begin{align}
&\phantom{+ }A^{ij}\eta _{ij}-\eta _{,i}g^{ij}F_{,u_{j}}-X(F)+\lambda F
+2A^{ij}\eta _{ui}u_{j}-A^{ij}\xi _{,ij}^{a}u_{a}-u_{i}\eta
_{u}g^{ij}F_{,u_{j}} +\xi _{,i}^{k}u_{k}g^{ij}F_{,u_{j}}  \notag \\
&+A^{ij}\eta _{uu}u_{i}u_{j}-2A^{ij}\xi
_{.,uj}^{k}u_{i}u_{k}+u_{i}u_{k}\xi_{,u}^{k}g^{ij}F_{,u_{j}} +A^{ij}\eta
_{u}u_{ij}-2A^{ij}\xi _{.,i}^{k}u_{jk}  \notag \\
&+(\xi ^{k}A_{,k}^{ij}+\eta A_{,u}^{ij})u_{ij}-\lambda A^{ij}u_{ij}
-A^{ij}\left( u_{ij}u_{a}+u_{i}u_{ja}+u_{ia}u_{j}\right) \xi
_{.,u}^{a}-u_{i}u_{j}u_{a}A^{ij}\xi _{uu}^{a}=0.  \label{Po.0}
\end{align}

We note that we cannot deduce the symmetry conditions before we select a
specific form for the function $F.$ However we may determine the conditions
which are due to the second derivative of $u$ because these terms do not
involve $F$. This observation significantly reduces the complexity of the
remaining symmetry condition. Following this observation we find the
condition
\begin{align*}
&A^{ij}\eta _{u}u_{ij}-A^{ij}(\xi _{.,i}^{k}u_{ja}+\xi
_{.,j}^{k}u_{ik})+(\xi ^{k}A_{,k}^{ij}+\eta A_{,u}^{ij})u_{ij}-\lambda
A^{ij}u_{ij} \\
& -A^{ij}\left( u_{ij}u_{a}+u_{i}u_{ja}+u_{ia}u_{j}\right) \xi
_{.,u}^{a}-u_{i}u_{j}u_{a}A^{ij}\xi _{uu}^{a}=0
\end{align*}%
which implies the equations
\begin{align*}
A^{ij}\left( u_{ij}u_{k}+u_{jk}u_{i}+u_{ik}u_{j}\right) \xi _{.,u}^{k}& =0,
\\
A^{ij}\eta _{u}u_{ij}-A^{ij}(\xi _{.,i}^{k}u_{jk}+\xi _{.,j}^{k}u_{ik})+(\xi
^{k}A_{,k}^{ij}+\eta A_{,u}^{ij})u_{ij}-\lambda A^{ij}u_{ij}& =0, \\
A^{ij}\xi _{uu}^{a}& =0.
\end{align*}%
The first equation is written as
\begin{equation}
A^{ij}\xi _{.,u}^{k}+A^{kj}\xi _{.,u}^{i}+A^{ik}\xi
_{.,u}^{j}=0\Leftrightarrow A^{(ij}\xi _{.,u}^{k)}=0.  \label{Po.1}
\end{equation}%
The second equation gives
\begin{equation}
A^{ij}\eta _{u}+\eta A_{,u}^{ij}+\xi ^{k}A_{,k}^{ij}-A^{kj}\xi
_{.,k}^{i}-A^{ik}\xi _{.,k}^{j}-\lambda A^{ij}=0.  \label{Po.2}
\end{equation}%
Therefore the last equation results in
\begin{equation}
\xi _{,uu}^{k}=0.  \label{Po.2a}
\end{equation}

It is straightforward to show that condition (\ref{Po.1}) implies $\xi
_{.,u}^{k}=0 $ which is a well known result.

From the analysis so far we obtain that for all second-order PDEs\ of the
form $A^{ij}u_{ij}-F(x^{i},u,u_{i})=0,$ for which at least one of the $%
A^{ij} $ is nonzero, the coefficients $\xi_{.,u}^{i}=0$ or $\xi^{i}=\xi
^{i}(x^{j}).$ Furthermore condition~(\ref{Po.2a}) is identically satisfied.
We consider that $A^{ij}~$is non-degenerate, furthermore the third symmetry
condition (\ref{Po.2}) can be written as follows%
\begin{equation}
L_{\xi ^{i}\partial _{i}}A^{ij}=\lambda A^{ij}-(\eta A^{ij})_{,u}.
\label{GPE.32}
\end{equation}%
This condition implies that for all second-order PDEs of the form $%
A^{ij}u_{ij}-F(x^{i},u,u_{i})=0$ for which $A^{ij}{}_{,u}=0$, i.e.~$%
A^{ij}=A^{ij}(x^{i}),$ the vector $\xi ^{i}\partial _{i}$ is a CKV of the
metric $A^{ij}$ with conformal factor ($\lambda -\eta _{u})(x).$

Moreover, using that $\xi _{,u}^{i}=0$ when at least one of the $A_{ij}\neq
0,$ the symmetry condition (\ref{Po.0}) is simplified as follows
\begin{align}
&\phantom{++ }A^{ij}\eta _{ij}-\eta _{,i}A^{ij}F_{,u_{j}}-X(F)+\lambda F
+2A^{ij}\eta _{ui}u_{j}-A^{ij}\xi _{,ij}^{a}u_{a}-u_{i}\eta
_{u}A^{ij}F_{,u_{j}}  \notag \\
&+\xi _{,i}^{k}u_{k}A^{ij}F_{,u_{j}} +A^{ij}\eta _{uu}u_{i}u_{j}+A^{ij}\eta
_{u}u_{ij}-2A^{ij}\xi _{.,i}^{k}u_{jk} +\big(\xi ^{k}A_{,k}^{ij}+\eta
A_{,u}^{ij}\big)u_{ij}-\lambda A^{ij}u_{ij}=0,  \label{GPE.30}
\end{align}%
which together with the condition (\ref{GPE.32}) are the complete set of
conditions \emph{for all} second-order PDEs of the form $%
A^{ij}u_{ij}-F(x^{i},u,u_{i})=0$ for which at least one of the $A_{ij}\neq 0$%
. This class of PDEs is quite general making the above result very useful.

In order to continue we need to assume that the function $F(x,u,u_{i})~$is
of a special form.

\subsection{The Lie point symmetry conditions for a linear function $%
F(x,u,u_{i})$}

\label{The Lie symmetry conditions for a linear function}

We consider the function $F(x,u,u_{i})$ to be of the form%
\begin{equation}
F(x,u,u_{i})=B^{k}(x,u)u_{k}+f(x,u),  \label{GPE.30a}
\end{equation}%
where $B^{k}(x,u)$ and $f(x,u)$ are arbitrary functions of their arguments.
In this case the PDE\ is of the form%
\begin{equation}
A^{ij}u_{ij}-B^{k}(x,u)u_{k}-f(x,u)=0.  \label{GPE.30.1}
\end{equation}

The Lie symmetries of this type of PDEs have been studied previously by
Ibragimov\cite{Ibrag}. Assuming that at least one of the $A_{ij}\neq 0$ the
Lie point symmetry conditions are (\ref{GPE.32}) and (\ref{GPE.30}).

Replacing $F(x,u,u_{1})$ in (\ref{GPE.30}) we obtain the following result
\cite{JGP1}:

\label{propPDE.1} The Lie point symmetry conditions for the second-order
PDEs of the form
\begin{equation}
A^{ij}u_{ij}-B^{k}(x,u)u_{k}-f(x,u)=0,  \label{GPE.40}
\end{equation}%
where at least one of the $A_{ij}\neq 0,$ are
\begin{gather}
A^{ij}(a_{ij}u+b_{ij})-(a_{,i}u+b_{,i})B^{i}-\xi
^{k}f_{,k}-auf_{,u}-bf_{,u}+\lambda f=0,  \label{GPE.42} \\
A^{ij}\xi _{,ij}^{k}-2A^{ik}a_{,i}+aB^{k}+auB_{,u}^{k}-\xi
_{,i}^{k}B^{i}+\xi ^{i}B_{,i}^{k}-\lambda B^{k}+bB_{,u}^{k}=0,
\label{GPE.43} \\
L_{\xi ^{i}\partial _{i}}A^{ij}=(\lambda -a)A^{ij}-\eta A^{ij}{}_{,u},
\label{GPE.44} \\
\eta =a(x^{i})u+b(x^{i}),  \label{GPE.45} \\
\xi _{,u}^{k} =0\Leftrightarrow \xi ^{k}(x^{i}).  \label{GPE.46}
\end{gather}%
From (\ref{GPE.44}) we infer that for all second-order PDEs of the form $%
A^{ij}u_{ij}-B^{k}(x,u)u_{k}-f(x,u)=0$ for which $A^{ij}{}_{,u}=0,$ the $\xi
^{i}(x^{j})$ is a CKV of the metric $A^{ij}.$ Also in this case for the
arbitrary function $\lambda $ holds $\lambda =\lambda (x^{i}).$ This result
establishes a relation between the Lie point symmetries of this type of PDEs
with the collineations of the metric defined by the coefficients $A_{ij}.$

Furthermore in case\footnote{%
The coordinates are $t,x^{i}$ where $i=1,\dots,n$.} when $%
A^{tt}=A^{tx^{i}}=0 $ and $A^{ij}$ is a non-degenerate metric we obtain that%
\begin{equation}
\xi _{,i}^{t}=0\Leftrightarrow \xi ^{t}(t).  \label{GPE.46a}
\end{equation}
These symmetry relations coincide with those given in \cite{Ibrag}. Finally
we note that equation (\ref{GPE.43}) can be written as
\begin{equation}
A^{ij}\xi _{,ij}^{k}-2A^{ik}a_{,i}+[\xi ,B]^{k}+(a-\lambda
)B^{k}+(au+b)B_{,u}^{k}=0.  \label{GPE.47}
\end{equation}%
Let us see now some applications of these results.

\section{Symmetries of the Poisson equation in a Riemannian space}

\label{PoissonSym}The Lie symmetries of the Poisson equation%
\begin{equation}
\Delta u-f\left( x^{i},u\right) =0,  \label{LE.01a}
\end{equation}%
where $\Delta =\frac{1}{\sqrt{\left\vert g\right\vert }}\frac{\partial }{%
\partial x^{i}}\left( \sqrt{\left\vert g\right\vert }g^{ij}\frac{\partial }{%
\partial x^{j}}\right)$ is the Laplace operator of the metric $g_{ij},$ for $%
f=f\left( u\right) $ have been given in \cite{Ibrag,Bozhkov}. Here we
generalize these results for the case $f=f\left( x^{i},u\right) $.

The Lie symmetry conditions (\ref{GPE.42})--(\ref{GPE.46}) for the Poisson
equation (\ref{LE.01a}) are
\begin{equation}
g^{ij}(a_{ij}u+b_{ij})-(a_{,i}u+b_{,i})\Gamma ^{i}-\xi
^{k}f_{,k}-auf_{,u}-bf_{,u}+\lambda f=0,  \label{WDW.01}
\end{equation}%
\begin{equation}
g^{ij}\xi _{,ij}^{k}-2g^{ik}a_{,i}+a\Gamma ^{k}-\xi _{,i}^{k}\Gamma ^{i}+\xi
^{i}\Gamma _{,i}^{k}-\lambda \Gamma ^{k}=0,  \label{WDW.02}
\end{equation}%
\begin{equation}
L_{\xi ^{i}\partial _{i}}g_{ij}=(a-\lambda )g_{ij},  \label{WDW.03}
\end{equation}%
\begin{equation}
\eta =a(x^{i})u+b(x^{i}),\quad\xi _{,u}^{k}=0.
\end{equation}%
Equation (\ref{WDW.02}) becomes (see \cite{Bozhkov})
\begin{equation}
g^{jk}L_{\xi }\Gamma _{.jk}^{i}=2g^{ik}a_{,i}.  \label{WDW.04}
\end{equation}%
From (\ref{WDW.03}), $\xi ^{i}$ is a CKV, then equations (\ref{WDW.04}) give
\begin{equation*}
\frac{2-n}{2}(a-\lambda )^{,i}=2a^{,i}\rightarrow \left( a-\lambda \right)
^{i}=\frac{4}{2-n}a^{,i}.
\end{equation*}

We define
\begin{equation}
\psi =\frac{2}{2-n}a+a_{0},
\end{equation}%
where $\psi =\frac{1}{2}\left( a-\lambda \right) $ is the conformal factor
of $\xi ^{i}$, i.e. $L_{\xi }g_{ij}=2\psi g_{ij}$. Furthermore, we have%
\begin{eqnarray*}
\left( 2-n\right) \lambda ^{i} &=&\left( 2-n\right) a^{i}-4a^{i}, \\
\left( 2-n\right) \lambda ^{i} &=&-\left( n+2\right) a^{i}.
\end{eqnarray*}%
Finally, from (\ref{WDW.01}), we have the constraint%
\begin{equation}
g^{ij}a_{i;j}u+g^{ij}b_{;ij}-\xi ^{k}f_{,k}-auf_{,u}+\lambda f-bf_{,u}=0.
\end{equation}

For $n=2,$ it holds that $~g^{jk}L_{\xi }\Gamma _{.jk}^{i}=0;$ this means
that $a_{,i}=0\rightarrow a=a_{0}.$ From (\ref{WDW.03}), $\xi ^{i}$ is a CKV
with conformal factor
\begin{equation}
2\psi =\left( a_{0}-\lambda \right)
\end{equation}%
and~$\lambda =a_{0}-2\psi .~$Finally, from (\ref{WDW.01}), we have the
constraint%
\begin{equation}
g^{ij}b_{;ij}-\xi ^{k}f_{,k}-a_{0}uf_{,u}+\left( a_{0}-2\psi \right)
f-bf_{,u}=0.
\end{equation}

Hence for the Lie symmetries of the Poisson equation in a general Riemannian
space we have the following theorem.

\begin{theorem}
\label{Theor}The Lie symmetries of the Poisson equation \eqref{LE.01a} are
generated from the CKVs of the metric~$g_{ij}$ defining the Laplace
operator, as follows

a) for $n>2$, the Lie symmetry vector is%
\begin{equation}
X=\xi ^{i}\big(x^{k}\big) \partial _{i}+\left( \frac{2-n}{2}\psi \big( x^{k}%
\big) u+a_{0}u+b\big(x^{k}\big) \right) \partial _{u},
\end{equation}%
where $\xi ^{i}\big(x^{k}\big) $ is a CKV with conformal factor $\psi \big(%
x^{k}\big) $ and the following condition holds%
\begin{equation}
\frac{2-n}{2}\Delta \psi u+g^{ij}b_{i;j}-\xi ^{k}f_{,k}-\frac{2-n}{2}\psi
uf_{,u}-\frac{2+n}{2}\psi f-bf_{,u}=0.  \label{KG.Eq0}
\end{equation}

b) for $n=2$, the Lie symmetry vector is
\begin{equation}
X=\xi ^{i}\big(x^{k}\big) \partial _{i}+\big( a_{0}u+b\big( x^{k}\big) \big) %
\partial _{u},
\end{equation}%
where $\xi ^{i}\big(x^{k}\big) $ is a CKV with conformal factor $\psi \big(%
x^{k}\big) $ and the following condition holds%
\begin{equation*}
g^{ij}b_{;ij}-\xi ^{k}f_{,k}-a_{0}uf_{,u}+\left( a_{0}-2\psi \right)
f-bf_{,u}=0.
\end{equation*}
\end{theorem}

\section{Lie Symmetries of the Klein--Gordon equation and CKVs}

\label{PoissonSym1}

In the special case, where $f\left( x^{i},u\right) =-V\left( x^{i}\right) u$%
, the Poisson equation (\ref{LE.01a}) is reduced to the Klein--Gordon
equation
\begin{equation}
\Delta u+V\big( x^{k}\big) u=0.  \label{KG.Eq1}
\end{equation}%
Therefore from theorem \ref{Theor} for the Lie symmetries of (\ref{KG.Eq1})
we have:

\begin{theorem}
\label{KG}The Lie point symmetries of the Klein--Gordon equation %
\eqref{KG.Eq1} are generated from the elements of the conformal algebra of
the metric~$g_{ij}$ defining the Laplace operator~$\Delta$. Equation~%
\eqref{KG.Eq1} admits the Lie symmetry vector%
\begin{equation}
X=\xi ^{i}\big(x^{k}\big) \partial _{i}+\left( \frac{2-n}{2}\psi \big( x^{k}%
\big) u+a_{0}u+b\big(x^{k}\big) \right) \partial _{u},  \label{KGT.01}
\end{equation}%
where $\xi ^{i}\big(x^{k}\big) $ is a CKV for the metric $g_{ij}$ with
conformal factor $\psi \big(x^{k}\big) $, the potential $V\big(
x^{k}\big) $ satisfies the condition%
\begin{equation}
\xi ^{k}V_{,k}+2\psi V-\frac{2-n}{2}\Delta \psi =0,  \label{KGT.02}
\end{equation}%
and the function $b\big(x^{k}\big) $ is a solution of \eqref{KG.Eq1}.
\end{theorem}

Of interest are the cases where $V\big(x^{k}\big) =0$ and $V\big(
x^{k}\big) =\frac{n-2}{4\left( n-1\right) }R$ , where $R$ is the Ricci
scalar of the metric $g_{ij}$ which defines the operator $\Delta $. In these
cases, the Klein--Gordon equation~(\ref{KG.Eq1}) becomes the Laplace
equation
\begin{equation}
\Delta u=0  \label{PE.9}
\end{equation}%
and the conformal invariant Laplace equation
\begin{equation}
\bar{L}_{g}u=0 ,  \label{PE.10}
\end{equation}%
where $\bar{L}_{g}=\Delta +\frac{n-2}{4\left( n-1\right) }R$. Therefore,
from theorem \ref{KG} for the Lie symmetries of equations (\ref{PE.9}) and (%
\ref{PE.10}) we have the following result

\begin{theorem}
\label{Laplace} The generic form of the Lie point symmetry of Laplace
equation \eqref{PE.9} and of the conformal invariant Laplace equation %
\eqref{PE.10} is the vector field \eqref{KGT.01} where for the Laplace
equation~\eqref{PE.9}, in a $n$-dimensional space with $n>2$, the conformal
factor $\psi \big(x^{k}\big) $ of the CKV $\xi ^{i}\big(x^{k}\big) $ is a
solution of \eqref{PE.9}.
\end{theorem}

\section{Reduction of Laplace equation in certain Riemannian spaces}

\label{Reduction}

From theorem \ref{Laplace} we have that the Lie symmetries of Laplace
equation (\ref{PE.9}) in a Riemannian space are generated from the CKVs (not
necessarily proper) whose conformal factor satisfies Laplace equation. This
condition is satisfied trivially by the KVs ($\psi =0),$ the HV ($\psi
_{;i}=0)$ and the sp.CKVs ($\psi _{;ij}=0).$ Therefore these vectors (which
span a subalgebra of the conformal algebra) are among the Lie symmetries of
Laplace equation. Concerning the proper CKVs it is not necessary that their
conformal factor satisfies the Laplace equation, therefore they may not
produce Lie symmetries for Laplace equation.

Furthermore, the special forms of the metric of a Riemannian space which
admits a gradient KV/HV or a sp.CKV are well known in the literature.
Therefore it is possible to study the application of the Lie symmetries of
the Laplace equation in these Riemannian spaces. In the following, we study
the reduction of Laplace equation in these general classes of Riemannian
spaces and we also study the origin of type II hidden symmetries. We assume
that the dimension~$n$ of the space is $n>2$.

\subsection{Reduction with a gradient KV/HV}

\label{grKV}

We consider the $(1+n)$-dimensional metric $g_{ij}~$with line element%
\begin{equation}
ds^{2}=dr^{2}+r^{2K}h_{AB}dy^{A}dy^{B},\quad h_{AB}=h_{AB}\left(
y^{C}\right),  \label{WH.03}
\end{equation}%
where $h_{AB}~$~is the metric of the $n$-dimensional space and $%
A,B,C=1,\dots,n$. For a general functional form of $h_{AB}$ and when $K=0$,
the metric (\ref{WH.03}) admits the gradient\ KV $\partial _{r}$, however
when $K=1$, the latter admits the gradient HV $r\partial _{r}~$(see \cite%
{Tupper1989}).

For the space with line element (\ref{WH.03}) Laplace equation (\ref{PE.9})
takes the form%
\begin{equation}
u_{,rr}+K\frac{n}{r}u_{,r}+\frac{1}{r^{2K}}~_{h}\Delta u=0,  \label{WH.04}
\end{equation}%
where $_{h}\Delta u=h^{AB}\left( y^{C}\right) u_{,AB}-\Gamma ^{A}\left(
y^{C}\right) u_{A}$ is the Laplace operator with metric $h_{AB}$.

Laplace equation (\ref{WH.04}) admits extra Lie point symmetries when $K=0,1$%
. In particular when $K=0$, the extra Lie point symmetry is the gradient KV $%
X_{KV}=\partial _{r}+\mu u\partial _{u}$ and when $K=1$ the extra Lie point
symmetry is the gradient HV $X_{HV}=r\partial _{r}+\mu u\partial _{u}$. We
will study the reduction of equation (\ref{WH.04}) using the zero-order
invariants of the symmetries $X_{KV}$ and $X_{HV}$.

The zero-order invariants of $X_{KV}$ are $\left\{ y^{A},e^{-\mu r}u\right\}
$ and of $X_{HV}$ are $\left\{ y^{A},r^{-\mu }u\right\} $. Hence by
replacing in equation (\ref{WH.04}) we find the reduced equation%
\begin{equation}
_{h}\Delta w+\mu ^{2}w=0,  \label{WH.05A}
\end{equation}%
where%
\begin{equation*}
u\left( r,y^{A}\right) =%
\begin{Bmatrix}
e^{\mu r}w\left( y^{A}\right) ~\text{ when~}~K=0 \\
r^{\mu }w\left( y^{A}\right) \text{ when~ }K=1%
\end{Bmatrix}%
.
\end{equation*}

Equation (\ref{WH.05A}) is the linear Klein--Gordon equation in the space
with metric $h_{AB}$. According to theorem \ref{KG} the Lie point symmetries
of the reduced equation (\ref{WH.05A}) follow from the CKVs of the metric $%
h_{AB}$.

The relation of the conformal algebra of the $n$ metric $h_{AB}$ and of the $%
1+n$ metric (\ref{WH.03}) have been studied in \cite{TNA}$.$ In particular
the KVs of the metrics $g_{ij}$ and $h_{AB}$ are the same. Furthermore, for $%
K=0,~$the $1+n$ metric $g_{ij}$ admits a HV\ if and only if the $n$ metric
admits one and if $_{n}H^{A}$ is the HV of the $n$ metric then the HV\ of
the $1+n$ metric is given by the expression{\samepage
\begin{equation}
_{1+n}H^{\mu }=z\delta _{z}^{\mu }+_{n}H^{A}\delta _{A}^{\mu }, \quad \mu
=x,1,\dots,n.
\end{equation}%
However, for $K=1$, the HV of the metric $g_{ij}$ is independent from that
of $h_{AB}$.}

Finally, the metric (\ref{WH.03}) admits proper CKVs if and only if the $n$
metric $h_{AB}$ admits gradient CKVs. This is because (\ref{WH.03}) is
conformally related with the decomposable metric
\begin{equation}
ds^{2\prime }=d\bar{r}^{2}+h_{AB}\left( y^{C}\right) dy^{A}dy^{B}.
\end{equation}%
which admits CKVs if and only if the $h_{AB}$ metric admits gradient CKVs.

The last implies, that Type II hidden symmetries are generated from the
elements of the (proper) conformal algebra of the $n$-dimensional metric $%
h_{AB}$ (for $K=0,1$) whose conformal factor is a solution of the
Klein--Gordon equation (\ref{WH.05A}), according to theorem~\ref{KG}.
Furthermore when $K=1$, the HV is a Type~II hidden symmetry.

In the following we will study the origin of type II hidden symmetries in
Riemannian spaces which admit a sp.CKV.

\subsection{Reduction with a sp.CKV}

\label{spCKV}

It is known \cite{HallspCKV}, that if an $n=m+1$-dimensional ($n>2)$
Riemannian space admits a non null sp.CKVs then also admits a gradient HV
and as many gradient KVs as the number of sp.CKVs. In these spaces there
exists always a coordinate system in which the metric is written in the form
\begin{equation}
ds^{2}=-dz^{2}+dR^{2}+R^{2}f_{AB}\left( y^{C}\right) dy^{A}dy^{B},
\label{LES.01}
\end{equation}%
where $~f_{AB}\left( y^{C}\right)$, $A,B,C,\ldots=1,2,\dots,m-1$ is an $(m-1)
$-dimensional metric. For a general metric $f_{AB}$ the $n$-dimensional
metric (\ref{LES.01}) admits a three-dimensional conformal algebra with
elements
\begin{equation*}
K_{G}=\partial _{z},\quad H=z\partial _{z}+R\partial _{R},\quad C_{S}=\frac{%
z^{2}+R^{2}}{2}\partial _{z}+zR\partial _{R},
\end{equation*}%
where $K_{G}$ is a gradient KV, $H$ is a gradient HV and $C_{S}$ is a sp.CKV
with conformal factor $\psi _{C_{S}}=z$. In these coordinates Laplace
equation (\ref{PE.9}) takes the form%
\begin{equation}
-u_{zz}+u_{RR}+\frac{( m-1) }{R}u_{R}+\frac{1}{R^{2}}~_{f}\Delta u_{AB}=0.
\label{LES.02}
\end{equation}

From theorem \ref{Laplace}, we have that the extra Lie point symmetries of (%
\ref{LES.02}) are the vectors%
\begin{equation*}
X^{1}=K_{G}+\mu _{G}X_{u},\quad X^{2}=H+\mu _{H}X_{u},\quad
X^{3}=C_{S}+2pzX_{u},
\end{equation*}%
where~$2p=\frac{1-m}{2}.$ The nonzero commutators of the extra Lie point
symmetries are%
\begin{equation*}
\left[ X^{1},X^{2}\right] =K_{G},\quad\left[ X^{2},X^{3}\right] =X^{3},\quad%
\left[ X^{1},X^{3}\right] =X^{2}+2pX_{u}.
\end{equation*}

The application of the Lie symmetries which are generated by the gradient KV
and the gradient HV have been studied in section \ref{grKV}. However we
would like to note that if we reduce the Laplace equation (\ref{LES.02}) by
use of the Lie symmetry $X^{1}$, the reduced equation admits the inherited
symmetry $X^{2}$ if and only if $\mu _{G}=0$. Furthermore the reduction with
the gradient HV leads to a PDE which does not admit inherited symmetries.
The resulting type II hidden symmetries follow from the results of section %
\ref{grKV}.

Before we reduce (\ref{LES.02}) with the symmetry generated by the sp.CKV~$%
X^{3}$, it is best to write the metric (\ref{LES.01}) in the coordinates $%
\left\{ x, R, y^{A}\right\}$, where the variable $x$ is defined by the
relation $z=\sqrt{R\left( R-x^{-1}\right)}$. In the new variables the Lie
symmetry $X^{3}$ becomes%
\begin{equation}
X^{3}=\sqrt{R\left( R-x^{-1}\right) }\left( R\partial _{R}+2pu\partial
_{u}\right) .
\end{equation}%
The zero-order invariants of $X^{3}$ in the new coordinates are $\left\{ x,
y^{A}, R^{-2p}u\right\} $. We choose $x,y^{A}$ to be the independent
variables and $w=w\left( x,y^{A}\right) $ to be the dependent one; that is,
the solution of the Laplace equation is in the form $u\left( x, R,
y^{A}\right) =R^{2p}w\left( x, y^{A}\right) $. Replacing in (\ref{LES.02})
we find the reduced equation
\begin{equation}
x^{2}w_{xx}+f^{AB}w_{AB}-\Gamma ^{A}w_{A}-2p\left( 2p+1\right) w=0.
\label{LES.07}
\end{equation}

For different values of the dimension $m$ of the metric $f_{AB}$, equation (%
\ref{LES.07}) can be written in the following forms%
\begin{eqnarray}
x^{2}w_{xx}+w_{yy}+\frac{1}{4}w &=&0,~\text{when }m=2,  \label{LES.07A} \\
_{\left( m=3\right) }\bar{\Delta}w &=&0,~\text{when }m=3,  \label{LES.07B} \\
_{\left( m\succeq 4\right) }\bar{\Delta}w-2p\left( 2p+1\right) V\left( \phi
\right) w &=&0,~\text{when }m\succeq 4,  \label{LES.07C}
\end{eqnarray}%
where $_{\left( m=3\right) }\bar{\Delta}$ is the Laplace operator for the
metric%
\begin{equation}
d\bar{s}_{\left( m=3\right) }^{2}=\frac{1}{x^{4}}dx^{2}+\frac{1}{x^{2}}%
f_{AB}dy^{A}dy^{B}.  \label{LES.07B1}
\end{equation}%
The Laplace operator $_{\left( m\succeq 4\right) }\bar{\Delta}$ is defined
by the metric%
\begin{equation}
d\bar{s}_{\left( m\succeq 4\right) }^{2}=d\phi ^{2}+\frac{1}{V}%
f_{AB}dy^{A}dy^{B},  \label{LES.07C1}
\end{equation}
where $V\left( \phi \right) =\dfrac{\left( 2-m\right) ^{2}}{\phi ^{2}}$ and $%
d\phi =\dfrac{1}{xV}dx.$

By applying the Lie symmetry condition (\ref{pr.04}) for equation (\ref%
{LES.07A}) we find that the generic Lie symmetry vector is
\begin{equation*}
X=\xi ^{x}\left( x,y\right) \partial _{x}+\xi ^{y}\left( x,y\right) \partial
_{y}+\left( a_{0}w+b\right) \partial _{w},
\end{equation*}%
where
\begin{eqnarray*}
\xi ^{x}\left( x,y\right) &=&c_{1}x+i\left( F_{1}\left( y+i\ln x\right)
+F_{2}\left( y-i\ln x\right) \right), \\
\xi ^{y}\left( x,y\right) &=&F_{2}\left( y-i\ln x\right) -F_{1}\left( y+i\ln
x\right).
\end{eqnarray*}%
$\xi ^{i}=\left( \xi ^{x},\xi ^{y}\right)$ is the generic CKV of the
two-dimensional metric $A^{ij}=\mathop{diag}\left( x^{2}, 1\right) $.
Therefore all the proper CKVs of the two-dimensional metric $A^{ij}$
generate type II hidden symmetries. Recall that the conformal algebra of a
two-dimensional space is infinite-dimensional. The Lie point symmetry $%
x\partial _{x}$ is the inherited symmetry $H$.

Furthermore for equations (\ref{LES.07B}) and (\ref{LES.07C}) from theorems~%
\ref{Laplace} and \ref{KG} we have that the type II hidden symmetries are
generated by the proper CKVs of the metrics (\ref{LES.07B1}) and (\ref%
{LES.07C1}), respectively, with conformal factors such as the condition (\ref%
{KGT.01}) holds. Finally equations (\ref{LES.07B}) and (\ref{LES.07C}) admit
the inherited Lie point symmetry $H$.

\section{Application of the reduction of the Laplace equation in Riemannian
spaces}

\label{Application}

In section \ref{Reduction}, we studied the reduction of Laplace equation and
the origin of type II hidden symmetries in general Riemannian spaces which
admit a gradient KV, a gradient HV and a sp.CKV. In this section we apply
these general results in order to study the reduction of Laplace equation
and the type II hidden symmetries in the case where the Laplace operator is
defined by (a)~the $n$-dimensional Minkowski spacetime $M^{n}$ (b)~the
four-dimensional conformally flat Bianchi~I spacetime which admits a
gradient KV.

\subsection{The $1+\left( n-1\right) $ wave equation}

We consider the $1+\left( n-1\right) $ wave equation
\begin{equation}
-u_{,tt}+u_{,zz}+\delta ^{AB}u_{,A}u_{,B}=0  \label{Apl.01}
\end{equation}%
which is the Laplace equation in the $n$-dimensional Minkowski spacetime $%
M^{n}~$ with$~n>3$.

The $n$-dimensional Minkowski spacetime admits a conformal group $G_{C}$ of
dimension $\dim G_{c}=\frac{1}{2}\left( n+1\right) \left( n+2\right) $.
Furthermore, the conformal group admits the following subalgebras:

a. $T^{n}$, the translation group of $n$ gradient KVs
\begin{equation*}
K_{G}^{1}=\partial _{t},\quad K_{G}^{z}=\partial _{z},\quad
K_{G}^{A}=\partial _{y^{A}};
\end{equation*}

b. $SO\left( n\right) ,$ the Lie group of $\frac{1}{2}n\left( n-1\right) $
rotations
\begin{equation*}
X_{R}^{1A}=Y^{\alpha }\partial _{t}+t\partial _{Y^{\alpha }},\quad
X_{R}^{\alpha \beta }=Y^{\beta }\partial _{Y^{\alpha }}-Y^{\alpha }\partial
_{Y^{\beta }},
\end{equation*}%
where $Y^{\alpha }=\left( z,y^{A}\right) $;

c. One gradient HV%
\begin{equation*}
H=t\partial _{t}+z\partial _{z}+y^{A}\partial _{y^{A}};
\end{equation*}

d. $G_{spC},$ the Lie group of $n$ sp. CKVs%
\begin{equation*}
X_{C}^{1}=\frac{1}{2}\left( t^{2}+\sum_{a}\left( Y^{\alpha }\right)
^{2}\right) \partial _{t}+t\sum_{a}\left( Y^{\alpha }\partial _{Y^{\alpha
}}\right),
\end{equation*}%
\begin{equation*}
X_{C}^{\alpha }=tY^{\alpha }\partial _{t}+\frac{1}{2}\left( Y^{\alpha
}+t^{2}-\sum_{\beta \neq \alpha }\left( Y^{\beta }\right) ^{2}\right)
\partial _{Y^{\alpha }}+Y^{\alpha }\sum_{\beta \neq \alpha }\left( Y^{\beta
}\partial _{Y^{\beta }}\right).
\end{equation*}

Therefore from theorem \ref{Laplace} we have that the Lie symmetries of the
wave equation (\ref{Apl.01}) are
\begin{equation*}
X_{u}=u\partial _{u},\quad K_{G}^{1},\quad K_{G}^{z},\quad K_{G}^{A},\quad
X_{R}^{1A},\quad X_{R}^{AB},\quad H,\quad X_{C}^{1}-tX_{u},\quad
X_{C}^{\alpha }-Y^{\alpha }X_{u}
\end{equation*}%
and that the nonzero commutators of the Lie symmetries are%
\begin{gather*}
\left[ K_{G}^{I},X_{R}^{IJ}\right] =-K_{G}^{J},\quad \left[ K_{G}^{I},H%
\right] =K_{G}^{I},\quad \left[ K_{G}^{I},X_{C}^{I}\right] =H-X_{u}, \\
\left[ K_{G}^{I},X_{C}^{J}\right] =X_{R}^{IJ},\quad \left[ H,X_{C}^{I}\right]
=X_{C}^{I},\quad \left[ X_{R}^{IJ},X_{C}^{I}\right] =X_{C}^{J},
\end{gather*}%
where $I=1,a$.

In order to apply the results of section \ref{Reduction} we reduce (\ref%
{Apl.01}) by the use of the Lie point symmetries generated by the gradient
KV $K_{G}^{z}$, the gradient HV $H$ and the special CKV $X_{C}^{1}$.

Using the zero-order invariants of the Lie symmetries $K_{G}^{z}+\mu X_{u}$
in (\ref{Apl.01}) the reduced equation is the linear Klein--Gordon equation
in the $M^{n-1}$ spacetime%
\begin{equation}
-w_{,tt}+\delta ^{AB}w_{,A}w_{,B}+\mu ^{2}w=0,  \label{Apl.02}
\end{equation}%
where $u\left( t,z,y^{A}\right) =e^{\mu z}w\left( t,y^{A}\right).$

For arbitrary constant $\mu,$ equation (\ref{Apl.02}) admits the Lie point
symmetries:

\begin{equation*}
K_{G}^{1},\quad K_{G}^{A},\quad X_{R}^{1A},\quad X_{R}^{AB},\quad H,\quad
w\partial _{w},
\end{equation*}%
which are the KVs and the HV of the $M^{n-1}$; these symmetries are
inherited symmetries. However when $\mu =0$, equation (\ref{Apl.02}) admits
the extra Lie point symmetries
\begin{equation*}
\bar{X}_{C}^{1}-\frac{1}{2}tw\partial _{w},\quad \bar{X}_{C}^{A}-\frac{1}{2}%
y^{A}w\partial _{w},
\end{equation*}%
where $\bar{X}_{C}^{1},~\bar{X}_{C}^{A}$ are the$~\left( n-1\right) $
sp.CKVs of $M^{n-1}$; these symmetries are type II hidden symmetries.

We continue with the reduction of equation (\ref{Apl.01}) by using the
invariants of the Lie point symmetry $H+\mu X_{u}$. The wave equation (\ref%
{Apl.01}) in hyperspherical coordinates $(r,\theta ^{A})$ is in the form of (%
\ref{WH.04}) with $K=1$ and $_{h}\Delta u$ is the Laplacian of the $n-1$
hyperbolic sphere~$S_{h}^{n-1}$ with line element%
\begin{equation}
ds_{h}^{2}=d\theta _{1}^{2}+\cosh ^{2}\theta _{1}\left( d\theta
_{2}^{2}+\cosh ^{2}\theta _{2}\left( d\theta _{3}^{2}+\cosh ^{2}\theta
_{3}(\dots) \right) \right) .  \label{Apl.03A}
\end{equation}

Hence, from section \ref{Reduction} we have that the reduced equation is the
linear Klein--Gordon equation in the hyperbolic sphere (\ref{Apl.03A})
\begin{equation}
_{h}\Delta w+\mu ^{2}w=0,  \label{Apl.03}
\end{equation}%
where $u\left( r,\theta ^{A}\right) =r^{\mu }w\left( \theta ^{A}\right) $.
The space $S_{h}^{n-1}$ is a maximally symmetric space and admits a~$\frac{1%
}{2}n(n+1)$-dimensional CKVs. In particular admits the $\frac{1}{2}n(n-1)$
KVs which form the $SO\left( n\right)$ Lie group and $n$ proper gradient
CKVs with conformal factor $\psi \left( \theta ^{A}\right) $ such as%
\begin{equation}
_{h}\Delta \psi =\frac{2-n}{4(n-1)}R_{h},
\end{equation}%
where $R_{h}$ is the Ricci scalar of the $S_{h}^{n-1}$ sphere which is a
constant. Therefore from theorem~\ref{Laplace} we have that for an arbitrary
constant $\mu ,$equation (\ref{Apl.03}) admits as Lie point symmetries the
elements of the $SO\left( n\right) $ and the linear symmetry $w\partial
_{w};~$these are inherited symmetries. However, when $\mu ^{2}=\frac{n-2}{%
4\left( n-1\right) }R_{h}$ the Klein--Gordon equation (\ref{Apl.03}) becomes
the conformal invariant Laplace equation and admits $\frac{1}{2}n(n+1)$ Lie
point symmetries, which span the conformal group of the space $S_{h}^{n-1}$.
Therefore reduction with the gradient HV\ leads to the conformal invariant
Laplace equation and the type II symmetries are generated by the proper CKVs
of hyperbolic sphere~$S_{h}^{n-1}$.

In order to reduce the wave equation (\ref{Apl.01}) with the invariants of
the Lie point symmetry $X_{C}^{1}-tX_{u}$ we use the coordinates~$\left(
t,R,\phi ^{K}\right) $ and equation (\ref{Apl.01}) becomes
\begin{equation}
u_{,.tt}-u_{,RR}-\frac{n-2}{R}u_{R}-\frac{1}{R^{2}}~_{h}\bar{\Delta}u=0,
\label{Apl.04}
\end{equation}%
where $_{h}\bar{\Delta}u$ is the Laplacian of the $(n-2)$-dimensional
hyperbolic sphere $S_{h}^{n-2}$. From section~\ref{Reduction} we have that
under the second coordinate transformation $t=\sqrt{R\left( R+\tau
^{-1}\right) }$ the reduced equation is written in the simpler form%
\begin{eqnarray}
_{\left( m=3\right) }\bar{\Delta}w &=&0,\quad \text{when }n=4,
\label{Apl.05} \\
_{\left( m\succeq 4\right) }\bar{\Delta}w-2p\left( 2p+1\right) \frac{\left(
2-m\right) ^{2}}{T^{2}}w &=&0,\quad \text{when }n>4,  \label{Apl.06}
\end{eqnarray}%
where~$u\left( \tau ,R,\phi ^{K}\right) =R^{2p}w\left( \tau ,y^{A}\right) $
and $_{\left( m=3\right) }\bar{\Delta}w$,~$_{\left( m\succeq 4\right) }\bar{%
\Delta}w~$are the Laplace operators for the $(n-1)$-dimensional flat space $%
ds_{\left( n-1\right) }^{2}$ with line element%
\begin{equation}
ds_{\left( n-1\right) }^{2}=dT^{2}-T^{2}ds_{h\left( n-2\right) }^{2}
\label{Apl.07}
\end{equation}%
for the cases where $n=4$ and $n>4$ respectively. $ds_{h\left( n-2\right)
}^{2}$ is the line element (\ref{Apl.03A}) of dimension $n-2$.

Therefore equation (\ref{Apl.05}) admits ten extra Lie symmetries which form
the conformal algebra of the flat space (\ref{Apl.07}), whereas equation (%
\ref{Apl.06}) admits as extra Lie point symmetries the gradient HV~$%
T\partial _{T}~$and the rotation group $SO\left( n\right) $ of the metric (%
\ref{Apl.07}) which is an inherited symmetry. Hence type II hidden
symmetries we have only when $n=4$ and they are generated by KVs the
translation group $T^{3}$ and the sp.CKVs of the flat space (\ref{Apl.07}).

The reduction of the $1+3$ wave equation has been done previously in \cite%
{AGA06} and our results coincide with theirs when $n=4$. However, our
approach is geometrical and we are able to recognize the type II hidden
symmetries from the form of the metric which defines the Laplace operator
without further calculations.

\subsection{Bianchi I (diagonal)\ spacetime}

The Bianchi I spacetime is defined as the four-dimensional spacetime which
admits as KVs the translation group $T^{3}$ of the 3D Euclidian space \cite%
{BRyan}. The line element of the Bianchi I spacetime in the invariant basis\
where the KVs are $\left\{ \partial _{x},\partial _{y},\partial _{z}\right\}
$ is:%
\begin{equation}
ds^{2}=-dt^{2}+A^{2}(t) dx^{2}+B^{2}(t) dy^{2}+C^{2}(t) dz^{2}.
\label{Apl.08}
\end{equation}
For general functions $A(t)$, $B(t)$, and $C(t)$ the KVs are non gradient.
We restrict our considerations to the case where $C(t)=1$ so that the KV $%
\partial_{z}$ is a gradient KV. In this case Laplace equation for the
Bianchi I spacetime (\ref{Apl.08}) becomes
\begin{equation}
-u_{,tt}+u_{,zz}+A^{-2}u_{,xx}+B^{-2}u_{,yy}-\left( \frac{\dot{A}}{A}+\frac{%
\dot{B}}{B}\right) u_{,t}=0.  \label{Apl.09}
\end{equation}

From theorem \ref{Laplace}, we have that the Laplace equation (\ref{Apl.09})
admits as extra Lie point symmetries the $T^{3}$ Lie group. We consider the
reduction of equation (\ref{Apl.09}) with the symmetry
\begin{equation*}
X_{I}=\partial _{z}+\mu u\partial _{u}.
\end{equation*}%
From the application of the zero order invariance of $X_{I}$ $\ $we have the
reduced equation%
\begin{equation}
_{h}\Delta w+\mu ^{2}w=0,  \label{Apl.10}
\end{equation}%
where $_{h}\Delta w$ is the Laplace operator for the three-dimensional
metric
\begin{equation}
ds_{\left( 3\right) }^{2}=-dt^{2}+A^{2}(t) dx^{2}+B^{2}( t) dy^{2}.
\label{Ap.11}
\end{equation}

The three-dimensional spacetime (\ref{Ap.11}) for general functions $A(t)$
and $B(t)$ admits the two-dimensional conformal algebra with elements the
KVs $\left\{ \partial _{x},\partial _{y}\right\} $ which are inherited Lie
point symmetries for the Klein--Gordon equation (\ref{Apl.10}).

Moreover, when $A^{2}\left( t\right) =B^{2}\left( t\right) $, the spacetime (%
\ref{Ap.11}) admits the KV $X_{I}^{3}=y\partial _{x}-x\partial _{y}$ which
is a Lie point symmetry of equation~(\ref{Apl.10}). However, $X_{I}^{3}$ is
an inherited symmetry because it is also a KV of (\ref{Apl.08}) and a Lie
point symmetry of Laplace equation~(\ref{Apl.09}). Furthermore, when $%
A^{2}(t) =B^{2}(t) =t^{2}$ the space $ds_{\left( 3\right) }^{2}$ admits the
HV (\ref{Ap.11}) which is a Lie point symmetry of equation (\ref{Apl.10})
when $\mu =0$. In this case $H$ is a inherited symmetry because in that case
metric (\ref{Apl.08}) admits the HV $H,$ which is also a Lie point symmetry
for (\ref{Apl.09}) with commutator $\left[ X_{I},H\right] =\partial _{z}.$

Finally, when $\left( A^{2}(t), B^{2}(t)\right) =\left( \sin ^{2}t,\cos
^{2}t\right) $ or $\left( A^{2}(t) ,B^{2}(t) \right) =\left( \sinh
^{2}t,\cosh ^{2}t\right) $ the three-dimensional spacetime (\ref{Ap.11}) is
a maximally symmetric spacetime, i.e.\ a space of constant non vanishing
curvature, and admits six KVs and four gradient CKVs (for the conformal
algebra of the Bianchi I spacetime see \cite{TsAp}). Hence equation (\ref%
{Apl.10}) for arbitrary $\mu $ admits six Lie point symmetries, the KVs.
However, when $\mu =-\frac{1}{8}, R_{\left( 3\right) }$ admits ten Lie point
symmetries which are generated by the conformal algebra of (\ref{Ap.11}),
where $R_{\left( 3\right)}$ is the Ricci scalar of the three-dimensional
space $ds_{\left( 3\right) }^{2}$. The type II hidden symmetries are the
symmetries which are generated by the gradient CKVs. The KVs are inherited
symmetries because they are also symmetries of equation (\ref{Apl.09}).

\section{Conclusion}

In this work we have studied the connection between the Lie point symmetries
of a general type of second-order PDEs and the collineations of the metric
defined by the PDE\ itself. In particular we wrote the Lie symmetry
conditions in a geometric form which include the Lie derivatives of
geometric objects. Therefore, we proved that the Lie point symmetries of the
PDEs of this family, which are linear in the first derivatives, are related
with the elements of the conformal algebra of the metric defined by the PDE
modulo a constraint condition depending again on the PDE. Important elements
in this family are the Poisson equation, the Klein--Gordon equation and the
Laplace equation.\ To these equations we have derived the Lie point
symmetries for various classes of Riemannian spaces. In particular we have
studied the type II\ hidden symmetries of the wave equation in Minkowski
space and in Bianchi I spacetimes. The methodology we followed is geometric
and can be applied to other types of PDEs as well as to this type of PDEs in
other Riemannian spaces.

\ack AP acknowledges financial support of INFN.

\end{document}